\documentclass[11pt]{article}
\usepackage{amsmath}
\usepackage{amssymb}
\usepackage{url}

\numberwithin{equation}{section}
\newtheorem{theorem}{Theorem}[section]
\newtheorem{Remark}[theorem]{Remark}
\newenvironment{remark}{\begin{Remark}\rm}{\end{Remark}}

\newcommand\dstyle\displaystyle
\newcommand\bLP{\\[\bigskipamount]}
\newcommand\sPP{\\[\smallskipamount]\indent}
\newcommand\RR{\mathbb{R}}
\newcommand\ZZ{\mathbb{Z}}
\newcommand\Znonneg{\ZZ_{\ge0}}
\newcommand\al\alpha
\newcommand\be\beta
\newcommand\de\delta
\newcommand\si\sigma
\newcommand\la\lambda
\newcommand\Ga{\Gamma}
\newcommand\De{\Delta}
\newcommand\La{\Lambda}
\newcommand\half{\frac12}
\newcommand\thalf{\tfrac12}
\newcommand\iy\infty
\newcommand\LHS{left-hand side}
\newcommand\RHS{right-hand side}
\newcommand{\hyp}[5]{\,\mbox{}_{#1}F_{#2}\!\left(
  \genfrac{}{}{0pt}{}{#3}{#4};#5\right)}

\begin{document}
\title{Generalizations of an integral for Legendre polynomials by
Persson and Strang}
\author{Enno Diekema\footnote{\tt e.diekema@gmail.com}\;
and
Tom H. Koornwinder\footnote{Korteweg-de Vries Institute,
University of Amsterdam, \tt T.H.Koornwinder@uva.nl}}
\date{}
\maketitle
\begin{abstract}
Persson and Strang (2003) evaluated the integral over $[-1,1]$ of
a squared odd degree Legendre polynomial divided by $x^2$ as being equal to 2.
We consider a similar integral for orthogonal polynomials with respect to
a general even orthogonality measure, with Gegenbauer and Hermite polynomials
as explicit special cases.
Next, after a quadratic transformation, we are led to the general
nonsymmetric case, with Jacobi and Laguerre polynomials as
explicit special cases. Examples of indefinite summation also occur
in this context. The paper concludes with a generalization of
the earlier results for Hahn polynomials. There some adaptations have to be
made in order to arrive at relatively nice explicit evaluations.
\end{abstract}
%
%
\section{Introduction}
The idea of this article comes from an integral
\begin{equation}
\int_{-1}^{1}\,\Big(\frac{P_{2n+1}(x)}{x}\Big)^2dx=2
\label{36}
\end{equation}
given by Persson and Strang \cite[(34)]{6}.
Here $P_n$ is a Legendre polynomial.
They prove the identity by deriving a first order recurrence for the
\LHS\ of \eqref{36}.
A natural question is if this integral can be generalized for other orthogonal polynomials.
We consider this problem first for orthogonal polynomials with respect
to an even orthogonality measure (the symmetric case),
a class which includes the Legendre polynomials, but also, for instance,
the Gegenbauer and Hermite polynomials.
The method of finding a first order recurrence for the integral still works,
but a method involving the Christoffel-Darboux formula (to some extent
equivalent to the earlier method) turns out to be more powerful in the
case of the Gegenbauer polynomials.

By applying a quadratic transformation to the symmetric case we arrive
at the idea of a further generalization in the case of a general
orthogonality measure. Kernel polynomials, defined in terms of the
Christoffel-Darboux kernel, enter the integral here. Explicit examples are
considered for Jacobi and Laguerre polynomials. Explicitly summable
indefinite sums naturally occur here as side results.

The original Persson-Strang integral was motivated by a particular
application.
Related motivations may also be given for our generalizations
discussed until now
(see also Remark \ref{62}).
But our last section on Hahn polynomials is
driven by the pure mathematical question how the Persson-Strang integral
and the related indefinite sum may generalize throughout the
Askey scheme. For the Hahn polynomials it turns out that the most
straightforward generalizations do not admit nice explicit evaluations.
We are able to make there an adaptation which admits
relatively nice evaluations
and which still has our results for Jacobi polynomials as a limit case.
\bLP
{\bf Acknowledgement}\quad
We thank Christian Krattenthaler and Michael Schlosser for helpful remarks.
\section{Preliminaries on orthogonal polynomials}
Let $\mu$ be a positive Borel measure on $\RR$ with infinite support
(or equivalently a nondecreasing function on $\RR$ with an infinite number
of points of increase) such that
$\int_\RR |x|^n\,d\mu(x)<\iy$ for all $n\in\Znonneg$. Let $p_n(x)$ ($n\in\Znonneg$) be a polynomial in $x$ of degree $n$ such that
\begin{equation}
\int_\RR p_m(x)\,p_n(x)\,d\mu(x)=h_n\,\de_{m,n}\qquad(m,n\in\Znonneg)
\label{1}
\end{equation}
for certain constants $h_n$ (necessarily positive).
The polynomials $p_n$ are called
{\em orthogonal polynomials} with respect to the measure $\mu$, see for instance \cite[\S5.2]{2}, \cite{3} or \cite{7}.
Up to constant nonzero factors they are uniquely determined by the above properties.
Let $k_n$ be the coefficient of $x^n$ in $p_n(x)$.
We assume that $p_0(x)=k_0=1$.

If $\mu$ has support within some closed interval $I$ then we say that
the $p_n$ are orthogonal polynomials with respect to $\mu$ on $I$.
In many examples we have $d\mu(x)=w(x)\,dx$ on $I$ with the
{\em weight function} $w$ a
nonnegative integrable function on $I$. Then \eqref{1} becomes:
\begin{equation*}
\int_I p_m(x)\,p_n(x)\,w(x)\,dx=h_n\,\de_{m,n}.
\label{42}
\end{equation*}
In many other examples $\mu$ is a discrete measure given by positive
weights $w_j$ on points $x_j$ ($j\in J$, a countably infinite set).
Then \eqref{36} becomes:
\begin{equation*}
\sum_{j\in J} p_m(x_j)\,p_n(x_j)\,w_j=h_n\,\de_{m,n}.
\label{43}
\end{equation*}
We may take $J$ finite, say $J=\{0,1,\ldots,N\}$. Then the $p_n$ are
well-defined by orthogonality for $n=0,1,\ldots,N$.

Orthogonal polynomials satisfy a three-term recurrence relation
(see \cite[Theorem 5.2.2]{2})
\begin{equation}
p_{n+1}(x)=(A_nx+B_n)p_n(x)-C_{n}p_{n-1}(x),\qquad p_{-1}(x)=0,\qquad p_0(x)=1.
\label{2}
\end{equation}
with $A_n,B_n,C_n$ real and $A_{n-1}A_nC_n>0$.
If the measure $\mu$ is {\em even} (i.e., invariant under reflection
with respect to 0) then $B_n=0$ for all $n$ in \eqref{2}.
If $\mu$ is even and $d\mu(x)=w(x)\,dx$ on $I$ then $I=-I$ and
$w(x)=w(-x)$ ($x\in I$).

It follows immediately from \eqref{2} that
\begin{equation}
A_n=\frac{k_{n+1}}{k_n},\quad
C_n=\frac{k_{n-1}h_nA_n}{k_nh_{n-1}}
=\frac{k_{n-1}k_{n+1}h_n}{k_n^2h_{n-1}}\,.
\label{23}
\end{equation}
Furthermore,
\begin{equation}
C_{2n-1}=-\,\frac{p_{2n}(0)}{p_{2n-2}(0)}\quad
\mbox{if $\mu$ is even.}
\label{24}
\end{equation}

We will also need the {\em Christoffel-Darboux formula}
(see \cite[Remark 5.2.2]{2})
\begin{equation}
\label{11}
K_n(x,y):=\sum_{k=0}^{n}\frac{p_k(x)p_k(y)}{h_k}=\frac{k_n}{k_{n+1}\,h_n}\,
\frac{p_{n+1}(x)p_n(y)-p_n(x)p_{n+1}(y)}{x-y}\quad(x\ne y).
\end{equation}
Note that $K_n$ is the kernel of the integral operator which projects
onto the space of polynomials of degree $\le n$. Thus
\begin{equation}
\int_\RR p(x)\,K_n(x,y)\,d\mu(x)=p(y)\qquad
\mbox{($p$ a polynomial of degree $\le n$).}
\label{57}
\end{equation}

If we let $y\to x$ in \eqref{11} then we obtain
\begin{equation}
\label{26}
K_n(x,x)=\sum_{k=0}^{n}\frac{p_k(x)^2}{h_k}=\frac{k_n}{k_{n+1}\,h_n}\,
\big(p_{n+1}'(x)p_n(x)-p_n'(x)p_{n+1}(x)\big).
\end{equation}
If the orthogonality measure $\mu$ is even then we have as a special case
of \eqref{26}:
\begin{equation}
K_{2n}(0,0)
=\sum_{k=0}^{n}\frac{p_{2k}(0)^2}{h_{2k}}=\frac{k_{2n}}{k_{2n+1}\,h_{2n}}\,
p_{2n+1}'(0)p_{2n}(0).
\label{27}
\end{equation}
\begin{remark}
If we are dealing with classical orthogonal polynomials $p_n$ then $p_n'$
is again a classical orthogonal polynomial, so \eqref{26} can be written
very explicitly. However, if we are dealing with orthogonal
polynomials $p_n$ in the Hahn class, i.e., where the polynomials
$\{\De p_n\}_{n\in\ZZ_{\ge1}}$ are again orthogonal (with
the difference operator $\De$ being defined by
$(\De f)(x):=f(x+1)-f(x)\,$), then we may better consider \eqref{11}
for $y=x-1$ (see \cite[formula after (2.11)]{9} for an analogous
observation in the $q$-case), which yields
\begin{equation}
K_n(x,x-1)=\sum_{k=0}^n\frac{p_k(x)\,p_k(x-1)}{h_k}=
\frac{k_n}{k_{n+1}\,h_n}\,\Big(p_n(x)\,(\De p_{n+1})(x-1)-
p_{n+1}(x)\,(\De p_n)(x-1)\Big).
\label{38}
\end{equation}
Formula \eqref{38} will take a very explicit form for polynomials of
Hahn class. If we consider \eqref{38} for Hahn polynomials then
a suitable limit case will yield formula \eqref{26} for Jacobi polynomials.
\end{remark}
\section{The Persson-Strang integral generalized: the symmetric case}
As a generalization of the integral of Persson and Strang we want to
compute the integral
\begin{equation}
I_n:=\int_{-a}^{a}\,\Big(\frac{p_{2n+1}(x)}{x}\Big)^2\,d\mu(x),
\label{28}
\end{equation}
where the $p_n$ are orthogonal polynomials with respect to
an even measure $\mu$ on the interval $[-a,a]$ or $(-\iy,\iy)$.
We will discuss two methods to solve this problem.
The first method, followed by Persson and Strang in the Legendre case,
is by a recurrence relation for $I_n$. The second method uses the
Christoffel-Darboux formula.

\subsection{The recurrence method}
\noindent Rewriting equation \eqref{2} with $B_n=0$ as
\begin{equation*}
\frac{p_{2n+1}(x)}{x}=A_{2n}p_{2n}(x)-C_{2n}\frac{p_{2n-1}(x)}{x}\,,
\end{equation*}
squaring, and integrating over the interval $[-a,a]$ with respect to $\mu$
gives:
\begin{multline*}
\int_{-a}^{a}\,\Big(\frac{p_{2n+1}(x)}{x}\Big)^2\,d\mu(x)
=A_{2n}^2\int_{-a}^{a}\,p_{2n}(x)^2\,d\mu(x)\\
-2A_{2n}C_{2n}\int_{-a}^{a}\,p_{2n}(x)\frac{p_{2n-1}(x)}{x}\,d\mu(x)
+C_{2n}^2\int_{-a}^{a}\,\Big(\frac{p_{2n-1}(x)}{x}\Big)^2\,d\mu(x)\\
=A_{2n}^2\,h_{2n}
+C_{2n}^2\int_{-a}^{a}\,\Big(\frac{p_{2n-1}(x)}{x}\Big)^2\,d\mu(x).
\end{multline*}
Hence we obtain the recurrence
\begin{equation}
I_n=C_{2n}^2 I_{n-1}+A_{2n}^2\,h_{2n}
\label{5}
\end{equation}
with starting value
\begin{equation}
I_0=\int_{-a}^{a}\,\Big(\frac{p_1(x)}{x}\Big)^2w(x)\,dx=k_1^2 h_0.
\label{7}
\end{equation}
The first order inhomogeneous linear recurrence relation \eqref{5}
with initial value \eqref{7} has a unique solution.
As observed by Persson and Strang, in the case of Legendre polynomials
the solution is $I_n=2$ for all $n$, since this solution satisfies
\eqref{5} and \eqref{7}.
\subsection{Using the Christoffel-Darboux formula}
\label{33}
Putting $y=0$ in the Christoffel-Darboux formula \eqref{11} gives:
\begin{equation}
K_{2n+1}(x,0)
=\sum_{k=0}^{2n+1}\frac{p_k(x)p_k(0)}{h_k}=K_{2n}(x,0)
=\sum_{k=0}^n\frac{p_{2k}(x)p_{2k}(0)}{h_{2k}}
=\frac{k_{2n}p_{2n}(0)}{k_{2n+1}h_{2n}}\,\frac{p_{2n+1}(x)}{x}\,.
\label{58}
\end{equation}
\begin{remark}
\label{62}
Combination of \eqref{58} with \eqref{57} yields
\begin{equation}
\frac{k_{2n}p_{2n}(0)}{k_{2n+1}h_{2n}}\,
\int_{-a}^a p(x)\,\frac{p_{2n+1}(x)}{x}\,d\mu(x)=p(0)\qquad
\mbox{($p$ a polynomial of degree $\le 2n+1$).}
\label{59}
\end{equation}
Thus the linear functional $\la\colon p\mapsto p(0)$ on the
finite Hilbert space of 
real-valued polynomials of degree $\le 2n+1$ with respect
to the inner product $\langle p,q\rangle:=\int_{-a}^a p(x)\,q(x)\,d\mu(x)$
gives $\la(p)$ as a constant times the inner product of
$p$ with the polynomial $x\mapsto p_{2n+1}(x)/x$. Therefore, the
square of the norm $\|\la\|^2$ of the linear functional $\la$ equals
a constant times the integral $I_n$ given by \eqref{28}.
This gives a motivation for trying to compute $I_n$ explicitly.

Another motivation considers the \LHS\ of \eqref{59} with $p$ being
a white noise signal on $[-a,a]$. Then this expression equals the
projection of $p$ on the subspace of polynomials of degree $\le 2n+1$,
evaluated at 0. This value is a random variable. The expectation of
the square of this value equals a constant times $I_n$. This is related
to the motivation in Persson \& Strang \cite[\S4]{6}.
\end{remark}

Now square the two sides of the last equality in \eqref{58}
and integrate over the orthogonality interval with respect to
$\mu$, where we use the orthogonality property.
As a result we obtain
\begin{equation}
\label{13}
I_n=\int_{-a}^{a}\Big(\frac{p_{2n+1}}{x}\Big)^2\,d\mu(x)
=\Big(\frac{k_{2n+1}h_{2n}}{k_{2n}p_{2n}(0)}\Big)^2
\sum_{k=0}^{n}\frac{p_{2k}(0)^2}{h_{2k}}\,.
\end{equation}
Then a very simple expression for $I_n$ can be obtained by substitution of
\eqref{27} in \eqref{13}:
\begin{equation}
\int_{-a}^{a}\Big(\frac{p_{2n+1}(x)}{x}\Big)^2\,d\mu(x)
=\frac{k_{2n+1}\,h_{2n}\,p_{2n+1}'(0)}{k_{2n}\,p_{2n}(0)}\,.
\label{34}
\end{equation}
\begin{remark}
The sum \eqref{13} is equivalent to the recurrence
\[
\Big(\frac{k_{2n}p_{2n}(0)}{k_{2n+1}h_{2n}}\Big)^2I_n
=\Big(\frac{k_{2n-2}p_{2n-2}(0)}{k_{2n-1}h_{2n-2}}\Big)^2I_{n-1}
+\frac{p_{2n}(0)^2}{h_{2n}}
\]
with starting value \eqref{7}.
The recurrence can be rewritten as
\begin{equation}
I_n=\Big(\frac{k_{2n-2}k_{2n+1}h_{2n}p_{2n-2}(0)}
{k_{2n-1}k_{2n}h_{2n-2}p_{2n}(0)}\Big)^2 I_{n-1}+
\Big(\frac{k_{2n+1}}{k_{2n}}\Big)^2 h_{2n}.
\label{25}
\end{equation}
In view of \eqref{23} and \eqref{24}, the recurrences \eqref{25} and
\eqref{5} are the same.
\end{remark}
\subsection{Example: Gegenbauer polynomials}
With the notation of \S\ref{33} let $\al>-1$ and take orthogonality measure
$d\mu(x):=(1-x^2)^\al\,dx$ on the interval $[-1,1]$. Then the $p_n$ are
{\em Gegenbauer polynomials} which we write as special Jacobi polynomials
(see \cite[10.9(4)]{10}):
\[
p_n(x)=P_n^{(\al,\al)}(x)=\frac{(\al+1)_n}{(2\al+1)_n}\,C_n^{\al+\half}(x).
\]
For the evaluation of the \RHS\ of \eqref{34} in this case we need
(see \cite[\S10.8, 10.9]{10}):
\[
k_n=\frac{(n+2\al+1)_n}{2^n n!}\,,\quad
h_n=\frac{2^{2\al+1}\,\Ga(n+\al+1)^2}{(2n+2\al+1)\,\Ga(n+2\al+1)\,n!}\,,
\]
\[
p_{2n}(0)=\frac{(\al+1)_{2n}}{(2\al+1)_{2n}}\,C_{2n}^{\al+\half}(0)
=(-1)^n\,\frac{(\al+n+1)_n}{2^{2n}n!}\,,
\]
\[
p_{2n+1}'(0)=(n+\al+1)\,P_{2n}^{(\al+1,\al+1)}(0)=
(-1)^n\,\frac{(\al+n+2)_n}{2^{2n}n!}\,.
\]
Then \eqref{34} yields
\begin{equation}
I_n=\int_{-1}^{1}\Big(\frac{P_{2n+1}^{(\al,\al)}(x)}{x}\Big)^2
(1-x^2)^\al\,dx
=\frac{2^{2\al+1}\Ga(2n+\al+2)^2}
{\Ga(2n+2\al+2)\,(2n+1)!}\,.
\label{35}
\end{equation}
As special cases of \eqref{35} we note:
\begin{itemize}
\item
$\al=0$, Legendre polynomials $P_n:=P_n^{(0,0)}$. Then $I_n=2$ and we recover
\eqref{36}.
\item
$\al=-\thalf$, Chebyshev polynomials of the first kind
$T_n:=\dstyle\frac{n!}{(1/2)_n}\,P_n^{(-1/2,-1/2)}$. Then
\[
I_n=\int_{-1}^{1}\Big(\frac{T_{2n+1}(x)}{x}\Big)^2
(1-x^2)^{-1/2}\,dx=(2n+1)\pi,
\]
This corresponds to \eqref{5} and \eqref{7} becoming
$I_n=I_{n-1}+2\pi$, $I_0=\pi$.
\item
$\al=\thalf$, Chebyshev polynomials of the second kind
$U_n:=\dstyle\frac{(n+1)!}{(3/2)_n}\,P_n^{(1/2,1/2)}$. Then
\[
I_n=\int_{-1}^{1}\Big(\frac{U_{2n+1}(x)}{x}\Big)^2
(1-x^2)^{1/2}\,dx=(2n+2)\pi,
\]
This corresponds to \eqref{5} and \eqref{7} becoming
$I_n=I_{n-1}+2\pi$, $I_0=2\pi$.
\end{itemize}
\subsection{Example: Hermite polynomials}
With the notation of \S\ref{33} take $d\mu(x):=e^{-x^2}\,dx$
on the interval $(-\iy,\iy)$. Then the $p_n$ are {\em Hermite polynomials},
\[
p_n(x)=H_n(x),
\]
for which we have (see \cite[\S10.13]{10}):
\[
h_n=\pi^{1/2}\,2^n\,n!\,,\quad
k_n=2^n,\quad
p_{2n}(0)=(-1)^n\,2^{2n}\,(\thalf)_n\,,\quad
p_{2n+1}'(0)=2(2n+1) p_{2n}(0).
\]
Then \eqref{34} yields
\begin{equation*}
I_n=\int_{-\iy}^{\iy}\Big(\frac{H_{2n+1}(x)}{x}\Big)^2
e^{-x^2}\,dx
=\pi^{1/2}\,2^{2n+2}\,(2n+1)!\,.
\end{equation*}

\section{Persson-Strang type integrals for general measures}
\subsection{Quadratic transformation}
Let the polynomials $p_n$ be orthogonal with respect to an even weight
function $w$ on the interval $(-a,a)$. For nonzero constants $C_n$
define polynomials $q_n$ of degree $n$ by
\begin{equation*}
p_{2n+1}(x)=c_n\,x\,q_n(x^2).
\end{equation*}
Then the polynomials $q_n$ are orthogonal with respect to the
measure $x^{1/2}\,w(x^{1/2})\,dx$ on the interval $[\,0,\sqrt a\,)$,
see Chihara \cite[Ch. I, \S8]{3}. We
can also rewrite the integral \eqref{28} in terms of the polynomials $q_n$:
\begin{equation*}
\int_{-a}^{a}\Big(\frac{p_{2n+1}(x)}{x}\Big)^2w(x)\,dx=
c_n^2 \int_{0}^{\sqrt a}q_n(x)^2\,x^{-1/2}\,w(x^{1/2})\,dx.
\end{equation*}
This suggests that, for an orthogonality measure $\mu$
on an interval $[0,a]$ and
for orthogonal polynomials $q_n$ with respect to the measure $x\,d\mu(x)$
on $[0,a]$ the integral
\[
\int_0^a q_n(x)^2\,d\mu(x)
\]
may have a nice evaluation. Moreover, we recognize the polynomials $q_n$
as {\em kernel polynomials} corresponding to the orthogonal polynomials
on $[0,a]$ with respect to measure $\mu$.
\subsection{Using kernel polynomials}
\label{37}
Let $\{p_n\}$ be a system of orthogonal polynomials with respect
to an orthogonality measure $\mu$ with support within $(-\iy,a]$. Let $k_n$ be the coefficient of $x^n$ in $p_n(x)$. Let $x_0\ge a$. For certain nonzero constants $c_n$ put
\begin{equation}
q_n(x):=c_n\,K_n(x_0,x)=
c_n\,\sum_{k=0}^n \frac{p_k(x_0)\,p_k(x)}{h_k}\,.
\label{3}
\end{equation}
Then $\{q_n\}$ is a system of orthogonal polynomials with respect to the orthogonality measure $(x_0-x)\,d\mu(x)$. These polynomials are called
{\em kernel polynomials}
(see \cite[\S5.6]{2}). Let $k_n'$ be the coefficient of $x^n$ in $q_n(x)$. Then
\[
c_n=\frac{k_n' h_n}{k_n p_n(x_0)}\,.
\]
Note as a special case of \eqref{3}:
\begin{equation}
c_n\,\sum_{k=0}^n \frac{p_k(x_0)^2}{h_k}=q_n(x_0).
\label{4}
\end{equation}
The kernel polynomial property of $q_n$ follows by combination
of \eqref{57} and \eqref{3}:
\begin{equation}
\int_\RR q_n(x)\,p(x)\,d\mu(x)=c_n\,p(x_0)\qquad
\mbox{($p$ a polynomial of degree $\le n$).}
\label{41}
\end{equation}
In particular, with $p=q_n$,
\begin{equation}
\int_\RR q_n(x)^2\,d\mu(x)=c_n\,q_n(x_0).
\label{9}
\end{equation}
If, for special choices of $\mu$ and $x_0$, we can explicitly evaluate $p_n(x_0)$, $q_n(x_0)$, $h_n$ and $c_n$, then \eqref{9} and \eqref{4} yield possibly interesting explicit evaluations of an integral and a finite sum, respectively.

Note that all formulas in this subsection remain unchanged if
$\mu$ has support within $[a,\iy)$ and if $x_0\le a$. Then
$\{q_n\}$ is a system of orthogonal polynomials with respect to the orthogonality measure $(x-x_0)\,d\mu(x)$. For $x_0$ in the interior
of the orthogonality interval for $\{p_n\}$, the orthogonality property
of $\{q_n\}$ persists, but then the orthogonality measure is no longer
positive.

More generally than \eqref{4} and \eqref{9} we will consider later in
this paper specializations for Hahn polynomials of the identities
\begin{equation}
\sum_{k=0}^n \frac{p_k(x_0)\,p_k(x_1)}{h_k}
=K_n(x_0,x_1)=\int_\RR K_n(x_0,x)\,K_n(x_1,x)\,d\mu(x).
\label{49}
\end{equation}
\subsection{Example: Jacobi polynomials}
Let $\al,\be>-1$ and take for the orthogonality measure $d\mu(x):=(1-x)^\al(1+x)^\be\,dx$ with support $[-1,1]$. Then the $p_n$ are
{\em Jacobi polynomials} (see \cite[\S10.8]{10})
\[
p_n(x)=P_n^{(\al,\be)}(x),
\]
which are normalized by their value at $x_0:=1$:
\[
p_n(1)=\frac{(\al+1)_n}{n!}\,.
\]
Then
\[
h_n=\frac{2^{\al+\be+1}}{2n+\al+\be+1}\,
\frac{\Ga(n+\al+1)\,\Ga(n+\be+1)}{\Ga(n+\al+\be+1)\,n!}\,,\qquad
k_n=\frac{(n+\al+\be+1)_n}{2^n\,n!}\,.
\]
Furthermore, the $q_n$ are Jacobi polynomials
\begin{equation*}
q_n(x)=P_n^{(\al+1,\be)}(x),
\end{equation*}
for which
\begin{equation*}
q_n(1)=\frac{(\al+2)_n}{n!}\,,\quad
k_n'=\frac{(n+\al+\be+2)_n}{2^n\,n!}\,,\quad{\rm hence}\quad
c_n=\frac{2^{\al+\be+1}\,\Ga(\al+1)\,\Ga(n+\be+1)}{\Ga(n+\al+\be+2)}\,.
\end{equation*}
Substitution of the expressions for $c_n$ and $q_n(1)$ in \eqref{9} yields:
\begin{equation}
\int_{-1}^1 P_n^{(\al+1,\be)}(x)^2\,
(1-x)^\al(1+x)^\be\,dx
=\frac{2^{\al+\be+1}}{\al+1}\,
\frac{\Ga(n+\al+2)\,\Ga(n+\be+1)}{\Ga(n+\al+\be+2)\,n!}\,.
\label{10}
\end{equation}
Formula \eqref{10} is given without proof in \cite[p.285, formula (6)]{5}.
After substitution of the explicit values of $p_k(1)$, $h_k$, $q_n(1)$ and
$c_n$, formula \eqref{4} can be written as:
\begin{equation}
\sum_{k=0}^n\frac{((\al+\be+3)/2)_k\,(\al+\be+1)_k\,(\al+1)_k}
{((\al+\be+1)/2)_k\,(\be+1)_k\,k!}
=\frac{(\al+2)_n\,(\al+\be+2)_n}{(\be+1)_n\,n!}\,.
\label{12}
\end{equation}
The \LHS\ of \eqref{12} can be written as a terminating very well poised hypergeometric series, by which \eqref{12} takes the form
\begin{equation}
\hyp54{\al+\be+1,1+\thalf(\al+\be+1),\al+1,n+\al+\be+2,-n}
{\thalf(\al+\be+1),\be+1,-n,n+\al+\be+2}1=
\frac{(\al+2)_n\,(\al+\be+2)_n}{(\be+1)_n\,n!}\,.
\label{15}
\end{equation}
Formula \eqref{15}, which we derived here from \eqref{4}, is also a special case of \cite[Corollary 3.4.3]{2} (which, in its turn is a terminating version of a degenerate case of Dougall's evaluation of a terminating
2-balanced very well poised ${}_7F_6(1)$, see (2.2.9) and (2.2.10) in \cite{2}).
\begin{remark}
\label{61}
The \LHS\ of \eqref{12} is an example of an
{\em indefinite sum}: a sum $\sum_{k=0}^n c_k$, where $c_k$ is a
{\em hypergeometric term} (i.e., $c_{k+1}/c_k$ is a rational function of $k$)
which does not depend on the upper limit $n$ of the sum. Moreover,
by \eqref{12} the sum can be evaluated for each $n$ as a hypergeometric term
$s_n$ (i.e., $s_{n+1}/s_n$ is a rational function of $n$).
In general, {\em Gosper's algorithm} can test whether an indefinite sum
of hypergeometric terms is summable with a hypergeometric term as sum,
and it explicitly gives this sum if it exists
(see \cite[Ch. 5]{14}).
See \eqref{45} for a more involved example of such an indefinite sum
(which has \eqref{12} as a limit case). Of course, as soon as we have
an explicit indefinite summation $\sum_{k=0}^n c_k=s_n$ then
an a posteriori proof can be immediately given by checking $c_0=s_0$ and
$s_n-s_{n-1}=c_n$.
\end{remark}
\subsection{Example: Laguerre polynomials}
Let $\al>-1$ and take for the orthogonality measure
$d\mu(x):=x^\al\,e^{-x}\,dx$ with support $[0,\iy)$.
Then the $p_n$ are {\em Laguerre polynomials}
(see \cite[\S10.12]{10})
\[
p_n(x)=L_n^\al(x),
\]
which are normalized by their value at $x_0:=0$:
\[
p_n(0)=\frac{(\al+1)_n}{n!}\,.
\]
Then
\[
h_n=\frac{\Ga(n+\al+1)}{n!}\,,\qquad
k_n=\frac{(-1)^n}{n!}\,.
\]
Furthermore, the $q_n$ are Laguerre polynomials
\begin{equation*}
q_n(x)=L_n^{\al+1}(x),
\end{equation*}
for which
\begin{equation*}
q_n(0)=\frac{(\al+2)_n}{n!}\,\quad
k_n'=\frac{(-1)^n}{n!}\,,\quad{\rm hence}\quad
c_n=\Ga(\al+1).
\end{equation*}
Substitution of the expressions for $c_n$ and $q_n(1)$ in \eqref{9}
yields:
\begin{equation}
\int_{0}^\iy L_n^{\al+1}(x)^2\,
x^\al\,e^{-x}\,dx
=\Ga(\al+1)\,\frac{(\al+2)_n}{n!}\,.
\label{31}
\end{equation}
A more general version of formula \eqref{31}, but still a specialization
of \eqref{41}, was earlier obtained by
Carlitz \cite[p.340]{8} (however with an erroneous factor $(-1)^n$ and
without a side condition that $m\le n$).
Yet earlier, Mayr \cite[\S3]{15} evaluated the integral
$\int_0^\iy e^{-\la x} L_r^a(\al x) L_s^b(\be x) x^{\si-1}\,dx$
as an Appell $F_2$ hypergeometric function. Specialization of his formula
puts the \LHS\ of \eqref{31} equal to
\begin{align*}
&\Ga(\al+1)\,\Big(\frac{(\al+2)_n}{n!}\Big)^2\,
F_2(\al+1,-n,-n,\al+2,\al+2;1,1)\\
&\quad=\Ga(\al+1)\,\Big(\frac{(\al+2)_n}{n!}\Big)^2\,
\sum_{m=0}^n\frac{(\al+1)_m\,(-n)_m}{m!\,(\al+2)_m}\,
\hyp21{\al+m+1,-n}{\al+2}1\\
&\quad=\Ga(\al+1)\,\Big(\frac{(\al+2)_n}{n!}\Big)^2\,
\sum_{m=0}^n\frac{(\al+1)_m\,(-n)_m}{m!\,(\al+2)_m}\,
\frac{(1-m)_n}{(\al+2)_n}\,,
\end{align*}
which equals the \RHS\ of \eqref{31}. In the last equality we used
the Chu-Vandermonde formula \cite[Corollary 2.2.3]{2}.

After substitution of the explicit values of $p_k(0)$, $h_k$, $q_n(0)$ and
$c_n$, formula \eqref{4} can be written as:
\begin{equation}
\sum_{k=0}^n\frac{(\al+1)_k}{k!}=\frac{(\al+2)_n}{n!}\,.
\label{32}
\end{equation}
The \LHS\ can be rewritten as the terminating hypergeometric series
\[
{}_2F_1(-n,\al+1;-n;1).
\]
Hence \eqref{32} is a special case of
the Chu-Vandermonde formula.
Of course, \eqref{32} can also be immediately checked
(see end of Remark \ref{61}).
\section{Further generalization in the case of Hahn polynomials}
In this section we use {\em Hahn polynomials}
\begin{equation}
p_n(x)=Q_n(x;\al,\be,N)=
\hyp32{-n,n+\al+\be+1,-x}{\al+1,-N}1\qquad
(n=0,1,\ldots,N)
\label{40}
\end{equation}
(see \cite[\S1.5]{12}) for a pilot study in order to see how
the general theory of \S\ref{37} can be made concrete for families
of orthogonal polynomials higher up in the
Askey scheme. We will take $x_0=N$, but then it will turn out that the
right-hand sides of \eqref{4} and \eqref{9} cannot be made explicit in
a simple form. Instead we will therefore consider \eqref{49}
with $x_0=N$, $x_1=N-1$.

Hahn polynomials satisfy the orthogonality relation
\[
\sum_{x=0}^N p_m(x)\,p_n(x)\,w_x=h_n\,\de_{m,n}\qquad
(m,n=0,1,\ldots,N)
\]
with
\[
w_x=\frac{(\al+1)_x}{x!}\,\frac{(\be+1)_{N-x}}{(N-x)!}
\]
and
\begin{equation*}
h_n=\frac{(\al+\be+2)_N}{N!}\,\frac{n+\al+\be+1}{2n+\al+\be+1}\,
\frac{(-1)^n\,n!}{(-N)_n}\,\frac{(\be+1)_n}{(\al+1)_n}\,
\frac{(N+\al+\be+2)_n}{(\al+\be+2)_n}\,.
\end{equation*}
We have also (with notation of \S\ref{37} and with $x_0=N$):
\begin{align*}
&k_n=\frac{(n+\al+\be+1)_n}{(\al+1)_n\,(-N)_n}\,,\quad
p_n(N)=\frac{(-1)^n\,(\be+1)_n}{(\al+1)_n}\,,\\
&q_n(x)=Q_n(x;\al,\be+1,N-1),\quad
k_n'=\frac{(n+\al+\be+2)_n}{(\al+1)_n\,(-N+1)_n}\,.
\end{align*}
Hence
\[
c_n=\frac{(\al+\be+2)_N}{N!}\,
\frac{(N+\al+\be+2)_n}{(\al+\be+2)_n}\,\frac{n!}{(-N+1)_n}\,.
\]
Then \eqref{3} takes the form
\begin{multline}
Q_n(x;\al,\be+1,N-1)=
\frac{(N+\al+\be+2)_n}{(\al+\be+2)_n}\,\frac{n!}{(-N+1)_n}\\
\times\sum_{k=0}^n\frac{2k+\al+\be+1}{k+\al+\be+1}\,
\frac{(\al+\be+2)_k}{(N+\al+\be+2)_k}\,
\frac{(-N)_k}{k!}\,Q_k(x;\al,\be,N).
\label{39}
\end{multline}

Instead of putting $x=N$ in \eqref{39} (like we obtained \eqref{4} from
\eqref{3}), we can better put $x=N-1$ in \eqref{39}.
Indeed, $Q_n(N;\al,\be+1,N-1)$ does not have a simple explicit expression,
but there is a simple expression
\[
q_n(N-1)=Q_n(N-1;\al,\be+1,N-1)
=\frac{(-1)^n\,(\be+2)_n}{(\al+1)_n}\,,
\]
while
\[
p_n(N-1)=Q_n(N-1;\al,\be,N)=\frac{(-1)^n\,(\be+1)_n}{(\al+1)_n}\,
\Big(1-\,\frac{n(n+\al+\be+1)}{(\be+1)N}\Big),
\]
as follows easily from \eqref{40}.
Thus we will specialize \eqref{49} for Hahn polynomials with $x_0=N$ and
$x_1=N-1$. Then, by recalling that $q_n(x)=c_n K_n(N,x)$ and by putting
\[
r_n(x):=c_n\,K_n(N-1,x),
\]
we can rewrite \eqref{49} as
\begin{equation}
c_n\,\sum_{k=0}^n\frac{p_k(N)\,p_k(N-1)}{h_k}=q_n(N-1)
=\frac1{c_n}\,\sum_{x=0}^N q_n(x)\,r_n(x)\,w_x.
\label{50}
\end{equation}
We will examine the two identities in \eqref{50} more closely in the next
two subsections.
\subsection{The first identity: terminating very well poised ${}_6F_5(-1)$}
The first identity in \eqref{50} can be rewritten more explicitly as follows.
\begin{multline}
\frac{(N+\al+\be+2)_n}{(\al+\be+2)_n}\,\frac{n!}{(-N+1)_n}\,
\sum_{k=0}^n\frac{2k+\al+\be+1}{k+\al+\be+1}\,
\frac{(\al+\be+2)_k}{(N+\al+\be+2)_k}\,
\frac{(-N)_k}{k!}\,
\frac{(-1)^k(\be+1)_k}{(\al+1)_k}\\
\times\Big(1-\,\frac{k(k+\al+\be+1)}{(\be+1)N}\Big)
=\frac{(-1)^n\,(\be+2)_n}{(\al+1)_n}\,.
\label{45}
\end{multline}
The \LHS\ of \eqref{45} can be written as a linear combination of
two terminating very well poised hypergeometric series of argument $-1$,
by which \eqref{45} takes the form
\begin{multline}
\hskip-0.5cm
\hyp65{\al+\be+1,1+\thalf(\al+\be+1),\be+1,-N,n+\al+\be+2,-n}
{\thalf(\al+\be+1),\al+1,N+\al+\be+2,-n,n+\al+\be+2}{-1}
-\frac{(\al+\be+2)(\al+\be+3)}{(N+\al+\be+2)(\al+1)}\\
\times
\hyp65{\al+\be+3,1+\thalf(\al+\be+3),\be+2,-N+1,n+\al+\be+3,-n+1}
{\thalf(\al+\be+3),\al+2,N+\al+\be+3,-n+1,n+\al+\be+3}{-1}\\
=\frac{(\be+2)_n\,(\al+\be+2)_n}{(\al+1)_n\,n!}\,
\frac{(-1)^n\,(-N+1)_n}{(N+\al+\be+2)_n}
\label{46}
\end{multline}
We can also write the \LHS\ of \eqref{46} as one single
terminating very well poised ${}_8F_7$ of argument $-1$:
\begin{multline}
\hyp87{\al+\be+1,1+\thalf(\al+\be+1),c+1,\al+\be+2-c,\be+1,-N,n+\al+\be+2,-n}
{\thalf(\al+\be+1),c,\al+\be+1-c,\al+1,N+\al+\be+2,-n,n+\al+\be+2}{-1}\\
=\frac{(\be+2)_n\,(\al+\be+2)_n}{(\al+1)_n\,n!}\,
\frac{(-1)^n\,(-N+1)_n}{(N+\al+\be+2)_n}\,,
\label{60}
\end{multline}
where
\begin{equation*}
c=\Big((\be+1)N+\tfrac14(\al+\be+1)^2\Big)^{1/2}.
\end{equation*}
Note that \eqref{15} (with $\al$ and $\be$ interchanged)
can be obtained as the limit of \eqref{46} or \eqref{60} for
$N\to\iy$.

Formula \eqref{46} can also be derived by combining a few identities
for hypergeometric functions given in the literature.
First apply Whipple's formula
\cite[(6.3)]{11},
\cite[Theorem 3.4.6]{2} to the two ${}_6F_5(-1)$'s in order to
rewrite the \LHS\ of \eqref{46} as a linear combination of two
0-balanced ${}_3F_2(1)$'s:
\begin{multline}
\frac{(-1)^n\,(\al+\be+2)_n}{n!}\left(
\hyp32{N+\al+1,n+\al+\be+2,-n}{N+\al+\be+2,\al+1}1\right.\\
\left.+\frac{n(n+\al+\be+2)}{(N+\al+\be+2)(\al+1)}\,
\hyp32{N+\al+1,n+\al+\be+3,-n+1}{N+\al+\be+3,\al+2}1\right).
\label{47}
\end{multline}
By the contiguous relation
(see \cite[(3.11)]{16})
\begin{equation}
\hyp32{a,b,c}{d,e}z=
\hyp32{a,b,c-1}{d,e}z+\frac{abz}{de}\,\hyp32{a+1,b+1,c}{d+1,e+1}z,
\label{51}
\end{equation}
which can be proved in a straightforward way,
expression \eqref{47} simplifies to
\begin{equation}
\frac{(-1)^n\,(\al+\be+2)_n}{n!}\,
\hyp32{N+\al,n+\al+\be+2,-n}{N+\al+\be+2,\al+1}1.
\label{48}
\end{equation}
By the Pfaff-Saalsch\"utz formula \cite[(2.2.8)]{2}
expression \eqref{48} is equal to the \RHS\ of~\eqref{46}.
\begin{remark}
Formula \eqref{51} can be extended to the more general contiguous relation
\begin{multline}
\hyp rs{a_1,\ldots,a_r}{b_1,\ldots,b_s}z
=\hyp rs{a_1,\ldots,a_{r-1},a_r-1}{b_1,\ldots,b_s}z\\
+\frac{a_1\ldots a_{r-1}z}{b_1\ldots b_s}\,
\hyp rs{a_1+1,\ldots,a_{r-1}+1,a_r}{b_1+1,\ldots,b_s+1}z.
\label{52}
\end{multline}
For the proof just write the summand for the power series expansion in $z$
of the \LHS\ as
\[
\frac{(a_1)_k\ldots(a_{r-1})_k(a_r)_{k-1}(a_r-1)}{(b_1)_k\ldots(b_s)_k\,k!}\,
z^k+
\frac{(a_1)_k\ldots(a_{r-1})_k(a_r)_{k-1}\,k}{(b_1)_k\ldots(b_s)_k\,k!}\,z^k.
\]
The $q$-analogue of
\eqref{52} is given by Krattenthaler \cite[(2.2)]{13}.
\end{remark}
\subsection{The second identity: the kernel polynomials $K_n(N-1,x)$}
In the second identity of \eqref{50} the only still unexplicit
expression is the polynomial $r_n(x)=K_n(N-1,x)$. As a kernel polynomial for the point $N-1$ it satisfies the property that
\begin{equation*}
\sum_{x=0}^N r_n(x)\,p(x)\,w_x=c_n\,p(N-1)\qquad
\mbox{($p$ a polynomial of degree $\le n$).}
\end{equation*}
For the evaluation in terms of Hahn polynomials of $r_n(x)$ we derive
the following:
\begin{align*}
r_n(x)&=c_n\sum_{k=0}^n \frac{p_k(N-1)\,p_k(x)}{h_k}
=\sum_{k=0}^n\Big(1-\frac{k(k+\al+\be+1)}{(\be+1)N}\Big)\,
\frac{p_k(N)\,p_k(x)}{h_k}\\
&=q_n(x)-\frac1{(\be+1)N}\,(\La q_n)(x),
\end{align*}
where (see \cite[(1.5.5)]{12})
\[
(\La f)(x):=(x+\al+1)(x-N)(\De f)(x)+x(x-\be-N-1)(\De f)(x-1).
\]
Now use (see \cite[(1.5.5), (1.5.7)]{12}) that
\begin{multline*}
(\La q_n)(x)=n(n+\al+\be+2)q_n(x)-(x+\al+1)(\De q_n)(x)\\
=n(n+\al+\be+2)Q_n(x;\al,\be+1,N-1)
+\frac{n(n+\al+\be+2)}{(\al+1)(N-1)}\,Q_{n-1}(x;\al+1,\be+2,N-2).
\end{multline*}
So we obtain that
\begin{multline}
r_n(x)=
\Big(1-\,\frac{n(n+\al+\be+2)}{(\be+1)\,N}\Big) Q_n(x;\al,\be+1,N-1)\\
-\,\frac{n(n+\al+\be+2)}{N(N-1)(\al+1)(\be+1)}\,(x+\al+1)\,
Q_{n-1}(x;\al+1,\be+2,N-2).
\end{multline}
Thus the second identity in \eqref{50} becomes:
\begin{multline}
\Big(1-\,\frac{n(n+\al+\be+2)}{(\be+1)\,N}\Big)
\sum_{x=0}^N \big(Q_n(x;\al,\be+1,N-1)\big)^2\,w_x
-\,\frac{n(n+\al+\be+2)}{N(N-1)(\al+1)(\be+1)}\\
\times
\sum_{x=0}^N Q_n(x;\al,\be+1,N-1)\,Q_{n-1}(x;\al+1,\be+2,N-2)\,(x+\al+1)\,w_x\\
=\frac{(\al+\be+2)_N}{N!}\,\frac{(-1)^n\,(\be+2)_n}{(\al+1)_n}\,
\frac{(N+\al+\be+2)_n}{(\al+\be+2)_n}\,\frac{n!}{(-N+1)_n}\,.
\label{53}
\end{multline}
\begin{remark}
In \eqref{53} replace $x$ by $NX$, divide both sides by $N^{\al+\be+1}$
and let $N\to\iy$. Then we obtain, at least formally, as a limit case
of \eqref{53} the identity
\begin{equation}
\int_{-1}^1 P_n^{(\al,\be+1)}(x)^2\,
(1-x)^\al(1+x)^\be\,dx
=\frac{2^{\al+\be+1}}{\be+1}\,
\frac{\Ga(n+\al+1)\,\Ga(n+\be+2)}{\Ga(n+\al+\be+2)\,n!}\,,
\label{54}
\end{equation}
which becomes \eqref{10} after an easy rewriting.
For the limit transition from \eqref{53} to \eqref{54} use that
\[
\lim_{N\to\iy} Q_n(Nx;\al,\be,N)=\frac{n!}{(\al+1)_n}\,
P_n^{(\al,\be)}(1-2x)
\]
(see \cite[(2.5.1)]{12}) and
\[
\lim_{k\to\iy}\frac{\Ga(k+a)}{\Ga(k+b)}\,k^{b-a}=1
\]
(see \cite[(1.4.3)]{2}).
Note that the second sum on the \LHS\ of \eqref{53} is killed in the
limit process.
\end{remark}

\quad\\
\begin{footnotesize}
E. Diekema, P.O.\ Box 2439, 7302 EP Apeldoorn, The Netherlands;
\sPP
email: {\tt e.diekema@gmail.com}
\bLP
T. H. Koornwinder, Korteweg-de Vries Institute, University of
Amsterdam,\sPP
P.O.\ Box 94248, 1090 GE Amsterdam, The Netherlands;
\sPP
email: {\tt T.H.Koornwinder@uva.nl}
\end{footnotesize}

\end{document}